\documentclass[reqno,12pt]{amsart}
\usepackage{amsmath,amsthm}
\usepackage[T1]{fontenc}
\usepackage{enumerate}
\theoremstyle{plain}

\newtheorem{expl}{Example}

\newcommand{\ra}{\rightarrow}

\begin{document}
\title{On projective systems of rational difference equations}
\author[Frank J. Palladino]{Frank J. Palladino}
\address{Department of Mathematics, University of Rhode Island,Kingston, RI 02881-0816, USA;}
\email{frank@math.uri.edu}
\date{October 14, 2011}

\begin{abstract}
\noindent 
We discuss first order systems of rational difference equations which have the property that lines through the origin are mapped into lines through the origin.
We call such systems projective systems of rational difference equations and we show a useful change of variables which helps us to understand the behavior in these cases. We include several examples to demonstrate the utility of this change of variables.

\end{abstract}
\maketitle

\section{Introduction}
We study first order systems of rational difference equations which have the property that lines through the origin are mapped into lines through the origin.
We call such systems projective systems of rational difference equations. Consider the most general first order system of linear fractional rational difference equations.
$$x_{n+1,1}=\frac{\alpha_{1}+\sum^{k}_{i=1}\beta_{1,i}x_{n,i}}{A_{1} +\sum^{k}_{i=1}B_{1,i}x_{n,i}},\quad n\in\mathbb{N},$$
$$x_{n+1,2}=\frac{\alpha_{2}+\sum^{k}_{i=1}\beta_{2,i}x_{n,i}}{A_{2} +\sum^{k}_{i=1}B_{2,i}x_{n,i}},\quad n\in\mathbb{N},$$
$$\begin{array}{c} 
\vdots \\
\end{array}$$
$$x_{n+1,k}=\frac{\alpha_{k}+\sum^{k}_{i=1}\beta_{k,i}x_{n,i}}{A_{k} +\sum^{k}_{i=1}B_{k,i}x_{n,i}},\quad n\in\mathbb{N}.$$
It is customary to assume that initial conditions and parameters are nonnegative. We will also make this assumption. However, we note that if proper care is taken to avoid problems which may arise from divison by zero, then the changes of variables we present in this article will also work in the case of complex parameters and complex initial conditions. 
In the following three situations this system is projective. If $\alpha_{i}=0=A_{i}$ for all $1\leq i \leq k$ then the system is projective. If $\alpha_{i}=0$ for all $1\leq i \leq k$, $A_{i}=A_{j}$ for all $i,j\in \{1,\dots , k\}$, and $B_{\ell i}=B_{\ell j}$ for all $i,j,\ell \in \{1,\dots , k\}$, then the system is projective.
If $A_{i}=0$ for all $1\leq i \leq k$, $\alpha_{i}=\alpha_{j}$ for all $i,j\in \{1,\dots , k\}$, and $\beta_{\ell i}=\beta_{\ell j}$ for all $i,j,\ell\in \{1,\dots , k\}$, then the system is projective.\par
In words, the three scenarios can be described as follows. The first scenario is when all of the constants in the numerator and denominator, $\alpha_{i}$ and $A_{i}$, are zero. Such systems of rational difference equations are commonly referred to as homogeneous systems of rational difference equations. Homogeneous systems of two rational difference equations have been studied extensively in \cite{huangknopf}.
The strategy employed in \cite{huangknopf} was to use the projective change of variables to carry out a reduction in the number of systems, from two to one. This allowed the authors of \cite{huangknopf} to determine the behavior in many cases of rational systems in the plane.
We will call such systems projective rational systems of homogeneous type.\par 
The second scenario is when all of the constants in the numerator, $\alpha_{i}$, are zero and the denominators are all the same. 
In this scenario, the projective change of variables is commonly used to embed the linear fractional system into a higher dimensional linear system. An example of this is the solution of the Riccati difference equation which can be found in \cite{riccati} and \cite{kulenovicladas}.
We will call such systems projective rational systems of linear type.\par 
The third scenario is when all of the constants in the denominator, $A_{i}$, are zero and the numerators are all the same. We will call such systems projective rational systems of hyperbolic type. We will focus on this scenario in several examples to show how these changes of variables can be used to simplify the problem considerably and determine global behavior. \par
In the next section of this article, we provide a change of variables for each of the three types of projective rational difference equation. We then, in the following section, give examples for each type of projective rational difference equation. We show how the change of variables can be used to determine the qualitative behavior of rational systems in several special cases.

\section{A change of variables for each scenario}
The general projective rational difference equation of homogeneous type is the following system of rational difference equations. We will assume positive initial conditions for our work here.
$$x_{n+1,1}=\frac{\sum^{k}_{i=1}\beta_{1,i}x_{n,i}}{\sum^{k}_{i=1}B_{1,i}x_{n,i}},\quad n\in\mathbb{N},$$
$$\begin{array}{c} 
\vdots \\
\end{array}$$
$$x_{n+1,k}=\frac{\sum^{k}_{i=1}\beta_{k,i}x_{n,i}}{\sum^{k}_{i=1}B_{k,i}x_{n,i}},\quad n\in\mathbb{N}.$$
Now we divide through by one of the variables, in this example we will use $x_{n+1,k}$ but this choice is arbitrary and any varaible may be used. In order to do this, we should have $k>1$. We get
$$\frac{x_{n+1,1}}{x_{n+1,k}}=\frac{\left(\sum^{k}_{i=1}\beta_{1,i}\frac{x_{n,i}}{x_{n,k}}\right)\left(\sum^{k}_{i=1}B_{k,i}\frac{x_{n,i}}{x_{n,k}}\right)}{\left(\sum^{k}_{i=1}\beta_{k,i}\frac{x_{n,i}}{x_{n,k}}\right)\left(\sum^{k}_{i=1}B_{1,i}\frac{x_{n,i}}{x_{n,k}}\right)},\quad n\in\mathbb{N},$$
$$\begin{array}{c} 
\vdots \\
\end{array}$$
$$\frac{x_{n+1,k-1}}{x_{n+1,k}}=\frac{\left(\sum^{k}_{i=1}\beta_{k-1,i}\frac{x_{n,i}}{x_{n,k}}\right)\left(\sum^{k}_{i=1}B_{k,i}\frac{x_{n,i}}{x_{n,k}}\right)}{\left(\sum^{k}_{i=1}\beta_{k,i}\frac{x_{n,i}}{x_{n,k}}\right)\left(\sum^{k}_{i=1}B_{k-1,i}\frac{x_{n,i}}{x_{n,k}}\right)},\quad n\in\mathbb{N}.$$
Using the change of variables $u_{n,i}=\frac{x_{n,i}}{x_{n,k}}$ we get the following.
$$u_{n+1,1}=\frac{\left(\beta_{1,k}+\sum^{k-1}_{i=1}\beta_{1,i}u_{n,i}\right)\left(B_{k,k}+\sum^{k-1}_{i=1}B_{k,i}u_{n,i}\right)}{\left(\beta_{k,k}+\sum^{k-1}_{i=1}\beta_{k,i}u_{n,i}\right)\left(B_{1,k}+\sum^{k-1}_{i=1}B_{1,i}u_{n,i}\right)},\quad n\in\mathbb{N},$$
$$\begin{array}{c} 
\vdots \\
\end{array}$$
$$u_{n+1,k-1}=\frac{\left(\beta_{k-1,k}+\sum^{k-1}_{i=1}\beta_{k-1,i}u_{n,i}\right)\left(B_{k,k}+\sum^{k-1}_{i=1}B_{k,i}u_{n,i}\right)}{\left(\beta_{k,k}+\sum^{k-1}_{i=1}\beta_{k,i}u_{n,i}\right)\left(B_{k-1,k}+\sum^{k-1}_{i=1}B_{k-1,i}u_{n,i}\right)},\quad n\in\mathbb{N}.$$
The general projective rational difference equation of linear type is the following system of rational difference equations. We will assume positive initial conditions for our work here.
$$x_{n+1,1}=\frac{\sum^{k}_{i=1}\beta_{1,i}x_{n,i}}{A_{1} +\sum^{k}_{i=1}B_{1,i}x_{n,i}},\quad n\in\mathbb{N},$$
$$\begin{array}{c} 
\vdots \\
\end{array}$$
$$x_{n+1,k}=\frac{\sum^{k}_{i=1}\beta_{k,i}x_{n,i}}{A_{1} +\sum^{k}_{i=1}B_{1,i}x_{n,i}},\quad n\in\mathbb{N}.$$
Now we divide through by one of the variables, in this example we will use $x_{n+1,k}$ but this choice is arbitrary and any varaible may be used. In order to do this, we should have $k>1$. We get
$$\frac{x_{n+1,1}}{x_{n+1,k}}=\frac{\sum^{k}_{i=1}\beta_{1,i}\frac{x_{n,i}}{x_{n,k}}}{\sum^{k}_{i=1}\beta_{k,i}\frac{x_{n,i}}{x_{n,k}}},\quad n\in\mathbb{N},$$
$$\begin{array}{c} 
\vdots \\
\end{array}$$
$$\frac{x_{n+1,k-1}}{x_{n+1,k}}=\frac{\sum^{k}_{i=1}\beta_{k-1,i}\frac{x_{n,i}}{x_{n,k}}}{\sum^{k}_{i=1}\beta_{k,i}\frac{x_{n,i}}{x_{n,k}}},\quad n\in\mathbb{N},$$
Using the change of variables $v_{n,i}=\frac{x_{n,i}}{x_{n,k}}$ we get the following.
$$v_{n+1,1}=\frac{\beta_{1,k}+\sum^{k-1}_{i=1}\beta_{1,i}v_{n,i}}{\beta_{k,k}+\sum^{k-1}_{i=1}\beta_{k,i}v_{n,i}},\quad n\in\mathbb{N},$$
$$\begin{array}{c} 
\vdots \\
\end{array}$$
$$v_{n+1,k-1}=\frac{\beta_{k-1,k}+\sum^{k-1}_{i=1}\beta_{k-1,i}v_{n,i}}{\beta_{k,k}+\sum^{k-1}_{i=1}\beta_{k,i}v_{n,i}},\quad n\in\mathbb{N}.$$
The general projective rational difference equation of hyperbolic type is the following system of rational difference equations. We will assume positive initial conditions for our work here.
$$x_{n+1,1}=\frac{\alpha_{1}+\sum^{k}_{i=1}\beta_{1,i}x_{n,i}}{\sum^{k}_{i=1}B_{1,i}x_{n,i}},\quad n\in\mathbb{N},$$
$$\begin{array}{c} 
\vdots \\
\end{array}$$
$$x_{n+1,k}=\frac{\alpha_{1}+\sum^{k}_{i=1}\beta_{1,i}x_{n,i}}{\sum^{k}_{i=1}B_{k,i}x_{n,i}},\quad n\in\mathbb{N}.$$
Now we divide through by one of the variables, in this example we will use $x_{n+1,k}$ but this choice is arbitrary and any varaible may be used. In order to do this, we should have $k>1$. We get
$$\frac{x_{n+1,1}}{x_{n+1,k}}=\frac{\sum^{k}_{i=1}B_{k,i}\frac{x_{n,i}}{x_{n,k}}}{\sum^{k}_{i=1}B_{1,i}\frac{x_{n,i}}{x_{n,k}}},\quad n\in\mathbb{N},$$
$$\begin{array}{c} 
\vdots \\
\end{array}$$
$$\frac{x_{n+1,k-1}}{x_{n+1,k}}=\frac{\sum^{k}_{i=1}B_{k,i}\frac{x_{n,i}}{x_{n,k}}}{\sum^{k}_{i=1}B_{k-1,i}\frac{x_{n,i}}{x_{n,k}}},\quad n\in\mathbb{N}.$$
Using the change of variables $w_{n,i}=\frac{x_{n,i}}{x_{n,k}}$ we get the following.
$$w_{n+1,1}=\frac{B_{k,k}+\sum^{k-1}_{i=1}B_{k,i}w_{n,i}}{B_{1,k}+\sum^{k-1}_{i=1}B_{1,i}w_{n,i}},\quad n\in\mathbb{N},$$
$$\begin{array}{c} 
\vdots \\
\end{array}$$
$$w_{n+1,k-1}=\frac{B_{k,k}+\sum^{k-1}_{i=1}B_{k,i}w_{n,i}}{B_{k-1,k}+\sum^{k-1}_{i=1}B_{k-1,i}w_{n,i}},\quad n\in\mathbb{N}.$$

\section{Examples and remarks}
We begin with a well known example which demonstrates the projective change of variables of linear type. Our first example provides a different change of variables for the Riccati equation than the standard change of variables found in \cite{kulenovicladas}.

\begin{expl}
Consider the difference equation
$$x_{n+1}=\frac{\alpha + \beta x_{n}}{A + B x_{n}},\;\; n=0,1,2, \dots,$$
with initial condition $x_{0} \in (0,\infty)$ and parameters $\alpha, \beta, A, B \in [0,\infty)$ such that $(A+B)(\alpha+\beta)\neq 0$. Inspired by the projective change of variables of linear type we may construct the following linear system in the plane,
$$y_{n+1}= \alpha z_{n} + \beta y_{n},\;\; n=0,1,2, \dots,$$
$$z_{n+1}= A z_{n} + B y_{n},\;\; n=0,1,2, \dots. $$
Choosing initial conditions $y_{0}=x_{0}$ and $z_{0}=1$ we have, $\frac{y_{n}}{z_{n}}=x_{n}$ for all $n\geq 0$. We prove this by induction, $n=0$ provides the base case. Notice that since $(A+B)(\alpha+\beta)\neq 0$ and $y_{0},z_{0}\in (0,\infty)$, $y_{n},z_{n}\in (0,\infty)$ for all $n\in\mathbb{N}$. 
Now, suppose $\frac{y_{n}}{z_{n}}=x_{n}$, then
$$x_{n+1}=\frac{\alpha + \beta x_{n}}{A + B x_{n}}=\frac{\alpha  + \beta \frac{y_{n}}{z_{n}}}{A  + B \frac{y_{n}}{z_{n}}}=\frac{\alpha z_{n} + \beta y_{n}}{A z_{n} + B y_{n}}=\frac{y_{n+1}}{z_{n+1}}.$$
So $\frac{y_{n}}{z_{n}}=x_{n}$ for all $n\geq 0$.
\end{expl}
In our next example, we use the projective change of variables of hyperbolic type and the results from \cite{merino} to determine the qualitative behavior of a system of 3 rational difference equations.
\begin{expl}
Consider the following system of three rational difference equations,
$$x_{n+1}=\frac{x_{n}}{Cy_{n}+Az_{n}},\;\; n=0,1,2, \dots,$$
$$y_{n+1}=\frac{x_{n}}{Dz_{n}},\;\; n=0,1,2, \dots,$$
$$z_{n+1}=\frac{x_{n}}{\beta x_{n}+\alpha z_{n}},\;\; n=0,1,2, \dots,$$
with initial conditions $x_{0},y_{0},z_{0}\in (0,\infty )$ and parameters $C,A,D,\beta ,\alpha \in (0,\infty )$. 
For this system, $\{x_{n}\}^{\infty}_{n=1}$ converges to
$$\bar{x}=\frac{2\left(\beta D - AD -\alpha\right)^{2}+ 4\alpha D \beta C + 2\left(\beta D - AD -\alpha\right)\sqrt{(\beta D - AD -\alpha)^{2}+4\alpha D \beta C}}{(2\beta C)\left( 2\alpha\beta C + \beta^{2} D - AD\beta -\alpha\beta + \beta\sqrt{(\beta D - AD -\alpha)^{2}+4\alpha D \beta C}\right)}.$$

$\{y_{n}\}^{\infty}_{n=1}$ converges to 
$$\bar{y}=\frac{\beta D - AD -\alpha + \sqrt{(\beta D - AD -\alpha)^{2}+4\alpha D \beta C}}{2\beta CD}.$$
  $\{z_{n}\}^{\infty}_{n=1}$ converges to
$$\bar{z}=\frac{\beta D - AD -\alpha + \sqrt{(\beta D - AD -\alpha)^{2}+4\alpha D \beta C}}{ 2\alpha\beta C + \beta^{2} D - AD\beta -\alpha\beta + \beta\sqrt{(\beta D - AD -\alpha)^{2}+4\alpha D \beta C}}.$$\par
\end{expl}

\begin{proof}

Dividing through by $z_{n+1}$ we have the following relationship,
$$\frac{x_{n+1}}{z_{n+1}}=\frac{\beta x_{n}+\alpha z_{n}}{Cy_{n}+Az_{n}}= \frac{\beta \frac{x_{n}}{z_{n}}+\alpha }{C\frac{y_{n}}{z_{n}}+A},\;\; n=0,1,2, \dots,$$
$$\frac{y_{n+1}}{z_{n+1}}=\frac{\beta x_{n}+\alpha z_{n}}{Dz_{n}}= \frac{\beta x_{n}}{Dz_{n}}+\frac{\alpha }{D},\;\; n=0,1,2, \dots.$$
Relabeling so that $u_{n}=\frac{x_{n}}{z_{n}}$ and $v_{n}=\frac{y_{n}}{z_{n}}$ we have,
$$u_{n+1}=\frac{\beta u_{n}+\alpha }{Cv_{n}+A},\;\; n=0,1,2, \dots,$$
$$v_{n+1}=\frac{\beta u_{n}}{D}+\frac{\alpha }{D},\;\; n=0,1,2, \dots.$$
So the sequence $\{u_{n}\}^{\infty}_{n=1}$ satisfies the following second order rational difference equation,
$$u_{n+1}=\frac{\beta u_{n}+\alpha }{\frac{\beta C u_{n-1}}{D}+A+\frac{\alpha }{D}},\;\; n=1,2, \dots,$$
with initial conditions $u_{0}=u_{0}$ and $u_{1}=\frac{\beta u_{0}+\alpha }{Cv_{0}+A}$. We know from \cite{merino} that $\{u_{n}\}^{\infty}_{n=1}$ converges to 
$$\bar{u}=\frac{\beta D - AD -\alpha + \sqrt{(\beta D - AD -\alpha)^{2}+4\alpha D \beta C}}{2\beta C}.$$
Thus, $\{v_{n}\}^{\infty}_{n=1}$ converges to 
$$\bar{v}=\frac{2\alpha C + \beta D - AD -\alpha + \sqrt{(\beta D - AD -\alpha)^{2}+4\alpha D \beta C}}{2CD}.$$
Thus, since $y_{n+1}=\frac{u_{n}}{D}$, $\{y_{n}\}^{\infty}_{n=1}$ converges to 
$$\bar{y}=\frac{\beta D - AD -\alpha + \sqrt{(\beta D - AD -\alpha)^{2}+4\alpha D \beta C}}{2\beta CD}.$$
Since $x_{n}=\frac{u_{n}y_{n}}{v_{n}}$, $\{x_{n}\}^{\infty}_{n=1}$ converges to
$$\bar{x}=\frac{2\left(\beta D - AD -\alpha\right)^{2}+ 4\alpha D \beta C + 2\left(\beta D - AD -\alpha\right)\sqrt{(\beta D - AD -\alpha)^{2}+4\alpha D \beta C}}{(2\beta C)\left( 2\alpha\beta C + \beta^{2} D - AD\beta -\alpha\beta + \beta\sqrt{(\beta D - AD -\alpha)^{2}+4\alpha D \beta C}\right)}.$$
Since $z_{n}=\frac{y_{n}}{v_{n}}$, $\{z_{n}\}^{\infty}_{n=1}$ converges to
$$\bar{z}=\frac{\beta D - AD -\alpha + \sqrt{(\beta D - AD -\alpha)^{2}+4\alpha D \beta C}}{ 2\alpha\beta C + \beta^{2} D - AD\beta -\alpha\beta + \beta\sqrt{(\beta D - AD -\alpha)^{2}+4\alpha D \beta C}}.$$\par
This example demonstrates the utility of the projective change of variables. If we were to try to analyze this system of three rational difference equations without using the projective change of variables and Merino's results in \cite{merino}, then we would be forced to overcome the same obstacles already overcome by Merino with the help of a lengthy Mathematica computation.
\end{proof}
In our next example, we use the projective change of variables of homogeneous type to reduce a 3 dimensional system which is neither competitive nor cooperative, to a cooperative system in the plane. We then use the available theory to determine the behavior for the reduced system. This allows us to determine the qualitative behavior of the original system of 3 rational difference equations.
\begin{expl}
Consider the following system of three rational difference equations,
$$x_{n+1}=\frac{x_{n}+y_{n}}{A_{1}z_{n}+x_{n}+y_{n}},\;\; n=0,1,2, \dots,$$
$$y_{n+1}=\frac{x_{n}+y_{n}}{A_{2}z_{n}+x_{n}+y_{n}},\;\; n=0,1,2, \dots,$$
$$z_{n+1}=\frac{\alpha z_{n}+x_{n}+y_{n}}{x_{n}+y_{n}},\;\; n=0,1,2, \dots,$$
with initial conditions $x_{0},y_{0},z_{0}\in (0,\infty )$ and positive parameters $\alpha$, $A_{1}$, and $A_{2}$. Define a polynomial, 
$$P(w)=-w^{3}+\left(2-A_{1}-A_{2}-\alpha\right)w^{2}+\left(A_{1}+A_{2}-\alpha A_{1}-\alpha A_{2}-A_{1}A_{2}\right)w-\alpha A_{1}A_{2}.$$
Further, define 
$$w_{m}=\frac{2\left(2-A_{1}-A_{2}-\alpha\right)+\sqrt{4\left(2-A_{1}-A_{2}-\alpha\right)^{2}+12\left(A_{1}+A_{2}-\alpha A_{1}-\alpha A_{2}-A_{1}A_{2}\right)}}{3}.$$
Whenever $P(w_{m})\leq 0$ or one of the following two conditions hold,
\begin{enumerate}
\item $A_{1}+A_{2}-\alpha A_{1}-\alpha A_{2}-A_{1}A_{2}< 0$ and 
$$2\left(2-A_{1}-A_{2}-\alpha\right)\leq \sqrt{-12\left(A_{1}+A_{2}-\alpha A_{1}-\alpha A_{2}-A_{1}A_{2}\right)},$$
\item $A_{1}+A_{2}-\alpha A_{1}-\alpha A_{2}-A_{1}A_{2}= 0$ and $A_{1}+A_{2}+\alpha \geq 2$,
\end{enumerate}
then 
$$lim_{n\rightarrow\infty} x_{n}=0,$$
$$lim_{n\rightarrow\infty} y_{n}=0,$$
$$lim_{n\rightarrow\infty} z_{n}=\infty.$$
Whenever $P(w_{m})>0$ and one of the following three conditions hold,
\begin{enumerate}
\item $A_{1}+A_{2}-\alpha A_{1}-\alpha A_{2}-A_{1}A_{2}< 0$ and 
$$2\left(2-A_{1}-A_{2}-\alpha\right)>\sqrt{-12\left(A_{1}+A_{2}-\alpha A_{1}-\alpha A_{2}-A_{1}A_{2}\right)},$$
\item $A_{1}+A_{2}-\alpha A_{1}-\alpha A_{2}-A_{1}A_{2}= 0$ and $A_{1}+A_{2}+\alpha < 2$,
\item $A_{1}+A_{2}-\alpha A_{1}-\alpha A_{2}-A_{1}A_{2}> 0$,
\end{enumerate}
then the equation
$$\overline{w}=\frac{\overline{w}^{2}}{\alpha+\overline{w}}\left(\frac{1}{A_{1}+\overline{w}}+\frac{1}{A_{2}+\overline{w}}\right),$$
 has three nonnegative solutions $\overline{w}=0,\overline{w}_{1},\overline{w}_{2}$ with $0<\overline{w}_{1}<\overline{w}_{2}$.\newline
If $\frac{x_{0}+y_{0}}{z_{0}}\in [0,\overline{w}_{1})$, then 
$$lim_{n\rightarrow\infty} x_{n}=0,$$
$$lim_{n\rightarrow\infty} y_{n}=0,$$
$$lim_{n\rightarrow\infty} z_{n}=\infty.$$
If $\frac{x_{0}+y_{0}}{z_{0}}=\overline{w}_{1}$,
then 
$$ x_{n}=\frac{\overline{w}_{1}}{A_{1}+\overline{w}_{1}},\;\;n\geq 1,$$
$$ y_{n}=\frac{\overline{w}_{1}}{A_{2}+\overline{w}_{1}},\;\;n\geq 1,$$
$$ z_{n}=1+\frac{\alpha }{\overline{w}_{1}},\;\;n\geq 1.$$
If $\frac{x_{0}+y_{0}}{z_{0}}\in (\overline{w}_{1},\infty)$, then 
$$lim_{n\rightarrow\infty} x_{n}=\frac{\overline{w}_{2}}{A_{1}+\overline{w}_{2}},$$
$$lim_{n\rightarrow\infty} y_{n}=\frac{\overline{w}_{2}}{A_{2}+\overline{w}_{2}},$$
$$lim_{n\rightarrow\infty} z_{n}=1+\frac{\alpha }{\overline{w}_{2}}.$$
\end{expl}

\begin{proof}
 Dividing through by $z_{n+1}$ we have the following relationship,
$$\frac{x_{n+1}}{z_{n+1}}=\frac{\left(\frac{x_{n}}{z_{n}}+\frac{y_{n}}{z_{n}}\right)^{2}}{\left(A_{1}+\frac{x_{n}}{z_{n}}+\frac{y_{n}}{z_{n}}\right)\left(\alpha+\frac{x_{n}}{z_{n}}+\frac{y_{n}}{z_{n}}\right)},\;\; n=0,1,2, \dots,$$
$$\frac{y_{n+1}}{z_{n+1}}=\frac{\left(\frac{x_{n}}{z_{n}}+\frac{y_{n}}{z_{n}}\right)^{2}}{\left(A_{2}+\frac{x_{n}}{z_{n}}+\frac{y_{n}}{z_{n}}\right)\left(\alpha+\frac{x_{n}}{z_{n}}+\frac{y_{n}}{z_{n}}\right)},\;\; n=0,1,2, \dots.$$
Relabeling so that $u_{n}=\frac{x_{n}}{z_{n}}$ and $v_{n}=\frac{y_{n}}{z_{n}}$ we have,
$$u_{n+1}=\frac{\left(u_{n}+v_{n}\right)^{2}}{\left(A_{1}+u_{n}+v_{n}\right)\left(\alpha+u_{n}+v_{n}\right)},\;\; n=0,1,2, \dots,$$
$$v_{n+1}=\frac{\left(u_{n}+v_{n}\right)^{2}}{\left(A_{2}+u_{n}+v_{n}\right)\left(\alpha+u_{n}+v_{n}\right)},\;\; n=0,1,2, \dots.$$
Let $w_{n}=u_{n}+v_{n}$, then
$$w_{n+1}=\frac{w_{n}^{2}}{\alpha+w_{n}}\left(\frac{1}{A_{1}+w_{n}}+\frac{1}{A_{2}+w_{n}}\right),\;\; n=0,1,2, \dots.$$
This is a first order difference equation in one variable where the function is monotone increasing, so we may use the stair step diagram to find the qualitative behavior of solutions. First, let us find the equilibria. We solve
$$\overline{w}=\frac{\overline{w}^{2}}{\alpha+\overline{w}}\left(\frac{1}{A_{1}+\overline{w}}+\frac{1}{A_{2}+\overline{w}}\right),$$
$$0=\overline{w}\left(\frac{\overline{w}\left(A_{1}+A_{2}+2\overline{w}\right)-\left(\alpha+\overline{w}\right)\left(A_{1}+\overline{w}\right)\left(A_{2}+\overline{w}\right)}{\left(\alpha+\overline{w}\right)\left(A_{1}+\overline{w}\right)\left(A_{2}+\overline{w}\right)}\right),$$
$$0=\overline{w}\left(\frac{-\overline{w}^{3}+\left(2-A_{1}-A_{2}-\alpha\right)\overline{w}^{2}+\left(A_{1}+A_{2}-\alpha A_{1}-\alpha A_{2}-A_{1}A_{2}\right)\overline{w}-\alpha A_{1}A_{2}}{\left(\alpha+\overline{w}\right)\left(A_{1}+\overline{w}\right)\left(A_{2}+\overline{w}\right)}\right).$$
So $\overline{w}=0$ is an equilibrium and by Descartes' rule of signs there are either 2 positive equilibira, or zero positive equilibria. Let us examine the polynomial,
$$P(w)=-w^{3}+\left(2-A_{1}-A_{2}-\alpha\right)w^{2}+\left(A_{1}+A_{2}-\alpha A_{1}-\alpha A_{2}-A_{1}A_{2}\right)w-\alpha A_{1}A_{2},$$
more closely. The derivative of this polynomial is
$$D(w)=-3w^{2}+2\left(2-A_{1}-A_{2}-\alpha\right)w+A_{1}+A_{2}-\alpha A_{1}-\alpha A_{2}-A_{1}A_{2}.$$
So, if $A_{1}+A_{2}-\alpha A_{1}-\alpha A_{2}-A_{1}A_{2}> 0$, then the root 
$$w_{m}=\frac{2\left(2-A_{1}-A_{2}-\alpha\right)+\sqrt{4\left(2-A_{1}-A_{2}-\alpha\right)^{2}+12\left(A_{1}+A_{2}-\alpha A_{1}-\alpha A_{2}-A_{1}A_{2}\right)}}{3},$$
is the only positive root of $D(w)$. Since $D(0)> 0$ and $lim_{w\ra\infty}D(w)=-\infty$ in this case, $w_{m}$ is a local maximum for $P(w)$ and the largest value of $P(w)$ for positive $w$.
Since $P(0)<0$, we have two possibilities in the case where $A_{1}+A_{2}-\alpha A_{1}-\alpha A_{2}-A_{1}A_{2}> 0$. 
If $A_{1}+A_{2}-\alpha A_{1}-\alpha A_{2}-A_{1}A_{2}> 0$ and $P(w_{m})>0$, then there are three equilibria $\overline{w}=0,\overline{w}_{1},\overline{w}_{2}$ with $0<\overline{w}_{1}<\overline{w}_{2}$. In this case,
$$f(w)=\frac{w^{2}}{\alpha+w}\left(\frac{1}{A_{1}+w}+\frac{1}{A_{2}+w}\right)-w$$
is negative on the intervals $(0,\overline{w}_{1})$ and $(\overline{w}_{2},\infty)$ and positive on the interval $(\overline{w}_{1},\overline{w}_{2})$. So, we have that the basin of attraction of $0$ is $[0,\overline{w}_{1})$, the basin of attraction of $\overline{w}_{1}$ is $\{\overline{w}_{1}\}$, and the basin of attraction of $\overline{w}_{2}$ is $(\overline{w}_{1},\infty)$.\par
If $A_{1}+A_{2}-\alpha A_{1}-\alpha A_{2}-A_{1}A_{2}> 0$ and $P(w_{m})=0$, then this contradicts Descartes' rule of signs since there is a single positive equilibrium in this case.
If $A_{1}+A_{2}-\alpha A_{1}-\alpha A_{2}-A_{1}A_{2}> 0$ and $P(w_{m})<0$, then the function
$$f(w)=\frac{w^{2}}{\alpha+w}\left(\frac{1}{A_{1}+w}+\frac{1}{A_{2}+w}\right)-w$$
is negative on $(0,\infty)$ so there is only one non-negative equilibrium $\overline{w}=0$ and this equilibrium is globally asymptotically stable.\par
If $A_{1}+A_{2}-\alpha A_{1}-\alpha A_{2}-A_{1}A_{2}= 0$ and $A_{1}+A_{2}+\alpha\geq 2$, then by Descartes' rule of signs $P(w)$ has no positive zeros, so $f(w)$ is negative on $(0,\infty)$ so there is only one non-negative equilibrium $\overline{w}=0$ and this equilibrium is globally asymptotically stable.\par
Now, whenever $A_{1}+A_{2}-\alpha A_{1}-\alpha A_{2}-A_{1}A_{2}= 0$ and $A_{1}+A_{2}+\alpha < 2$, then
$$D\left(\frac{2-A_{1}-A_{2}-\alpha}{3} \right)> 0,$$
 and $lim_{w\ra\infty}D(w)=-\infty$. Thus, in this case, $w_{m}$ is a local maximum for $P(w)$ and the largest value of $P(w)$ for positive $w$.
Since $P(0)<0$, we have two possibilities in the case where $A_{1}+A_{2}-\alpha A_{1}-\alpha A_{2}-A_{1}A_{2}= 0$ and $A_{1}+A_{2}+\alpha < 2$. 
If $A_{1}+A_{2}-\alpha A_{1}-\alpha A_{2}-A_{1}A_{2}= 0$, $A_{1}+A_{2}+\alpha < 2$, and $P(w_{m})>0$, then there are three equilibria $\overline{w}=0,\overline{w}_{1},\overline{w}_{2}$ with $0<\overline{w}_{1}<\overline{w}_{2}$. In this case,  
$f(w)$ is negative on the intervals $(0,\overline{w}_{1})$ and $(\overline{w}_{2},\infty)$ and positive on the interval $(\overline{w}_{1},\overline{w}_{2})$. So, the basin of attraction of $0$ is $[0,\overline{w}_{1})$, the basin of attraction of $\overline{w}_{1}$ is $\{\overline{w}_{1}\}$, and the basin of attraction of $\overline{w}_{2}$ is $(\overline{w}_{1},\infty)$.\par
If $A_{1}+A_{2}-\alpha A_{1}-\alpha A_{2}-A_{1}A_{2}= 0$, $A_{1}+A_{2}+\alpha < 2$, and $P(w_{m})=0$, then this contradicts Descartes' rule of signs since there is a single positive equilibrium in this case.
If $A_{1}+A_{2}-\alpha A_{1}-\alpha A_{2}-A_{1}A_{2}= 0$, $A_{1}+A_{2}+\alpha < 2$, and $P(w_{m})<0$, then the function $f(w)$
is negative on $(0,\infty)$ so there is only one non-negative equilibrium $\overline{w}=0$ and this equilibrium is globally asymptotically stable.\par
If $A_{1}+A_{2}-\alpha A_{1}-\alpha A_{2}-A_{1}A_{2}< 0$ and 
$$2\left(2-A_{1}-A_{2}-\alpha\right)<\sqrt{-12\left(A_{1}+A_{2}-\alpha A_{1}-\alpha A_{2}-A_{1}A_{2}\right)},$$
then $D(w)$ has no positive real roots. Since $D(0)$ is clearly negative, we get that $D(w)$ is negative on $(0,\infty)$. Thus, since $P(0)$ is negative, $P(w)$ is negative on $(0,\infty)$. Thus, $f(w)$ is negative on $(0,\infty)$. 
So there is only one non-negative equilibrium $\overline{w}=0$ and this equilibrium is globally asymptotically stable.\par
If $A_{1}+A_{2}-\alpha A_{1}-\alpha A_{2}-A_{1}A_{2}< 0$ and 
$$2\left(2-A_{1}-A_{2}-\alpha\right)=\sqrt{-12\left(A_{1}+A_{2}-\alpha A_{1}-\alpha A_{2}-A_{1}A_{2}\right)},$$
then $D(w)$ has only one positive real root. Since $D(0)$ is clearly negative and $lim_{w\rightarrow\infty}=-\infty$, we get that $D(w)$ is nonpositive on $(0,\infty)$. Thus, since $P(0)$ is negative, $P(w)$ is negative on $(0,\infty)$. Thus, $f(w)$ is negative on $(0,\infty)$. 
So there is only one non-negative equilibrium $\overline{w}=0$ and this equilibrium is globally asymptotically stable.\par
If $A_{1}+A_{2}-\alpha A_{1}-\alpha A_{2}-A_{1}A_{2}< 0$ and 
\begin{equation}
2\left(2-A_{1}-A_{2}-\alpha\right)>\sqrt{-12\left(A_{1}+A_{2}-\alpha A_{1}-\alpha A_{2}-A_{1}A_{2}\right)},
\end{equation}
then $D(w)$ has two positive real roots. $D(0)$ is clearly negative and 
$$D\left(\frac{2-A_{1}-A_{2}-\alpha}{3} \right)= \frac{(2-A_{1}-A_{2}-\alpha)^{2}}{3}-A_{1}+A_{2}-\alpha A_{1}-\alpha A_{2}-A_{1}A_{2}.$$
Thus by the inequality (1), $$D\left(\frac{2-A_{1}-A_{2}-\alpha}{3} \right)>0.$$
Moreover,
 $lim_{w\rightarrow\infty}=-\infty$. Thus, in this case, $w_{m}$ is a local maximum for $P(w)$ and either $P(0)$ or $P(w_{m})$ is the largest value of $P(w)$ for nonnegative $w$.
Since $P(0)<0$, we have two possibilities in the case where $A_{1}+A_{2}-\alpha A_{1}-\alpha A_{2}-A_{1}A_{2}< 0$ and (1) holds. 
If $A_{1}+A_{2}-\alpha A_{1}-\alpha A_{2}-A_{1}A_{2}< 0$, (1) holds, and $P(w_{m})>0$, then there are three equilibria $\overline{w}=0,\overline{w}_{1},\overline{w}_{2}$ with $0<\overline{w}_{1}<\overline{w}_{2}$. In this case,  
$f(w)$ is negative on the intervals $(0,\overline{w}_{1})$ and $(\overline{w}_{2},\infty)$ and positive on the interval $(\overline{w}_{1},\overline{w}_{2})$. So, the basin of attraction of $0$ is $[0,\overline{w}_{1})$, the basin of attraction of $\overline{w}_{1}$ is $\{\overline{w}_{1}\}$, and the basin of attraction of $\overline{w}_{2}$ is $(\overline{w}_{1},\infty)$.\par
If $A_{1}+A_{2}-\alpha A_{1}-\alpha A_{2}-A_{1}A_{2}< 0$, (1) holds, and $P(w_{m})=0$, then this contradicts Descartes' rule of signs since there is a single positive equilibrium in this case.
If $A_{1}+A_{2}-\alpha A_{1}-\alpha A_{2}-A_{1}A_{2}< 0$, (1) holds, and $P(w_{m})<0$, then the function $f(w)$
is negative on $(0,\infty)$ so there is only one non-negative equilibrium $\overline{w}=0$ and this equilibrium is globally asymptotically stable.\par
Thus summarizing the results obtained above via the stair step method from \cite{kulenovicladas}, we get the following result.
Whenever $P(w_{m})>0$ and one of the following three conditions hold,
\begin{enumerate}
\item $A_{1}+A_{2}-\alpha A_{1}-\alpha A_{2}-A_{1}A_{2}< 0$ and 
$$2\left(2-A_{1}-A_{2}-\alpha\right)>\sqrt{-12\left(A_{1}+A_{2}-\alpha A_{1}-\alpha A_{2}-A_{1}A_{2}\right)},$$
\item $A_{1}+A_{2}-\alpha A_{1}-\alpha A_{2}-A_{1}A_{2}= 0$ and $A_{1}+A_{2}+\alpha < 2$,
\item $A_{1}+A_{2}-\alpha A_{1}-\alpha A_{2}-A_{1}A_{2}> 0$,
\end{enumerate}
then there are three equilibria $\overline{w}=0,\overline{w}_{1},\overline{w}_{2}$ with $0<\overline{w}_{1}<\overline{w}_{2}$. The basin of attraction of $0$ is $[0,\overline{w}_{1})$, the basin of attraction of $\overline{w}_{1}$ is $\{\overline{w}_{1}\}$, and the basin of attraction of $\overline{w}_{2}$ is $(\overline{w}_{1},\infty)$.\par
Since $w_{n}=\frac{x_{n}+y_{n}}{z_{n}}$, we have:
$$x_{n+1}=\frac{1}{\frac{A_{1}}{w_{n}}+1},\;\; n=0,1,2, \dots,$$
$$y_{n+1}=\frac{1}{\frac{A_{2}}{w_{n}}+1},\;\; n=0,1,2, \dots,$$
$$z_{n+1}=1+\frac{\alpha }{w_{n}},\;\; n=0,1,2, \dots.$$
Thus in this case if $\frac{x_{0}+y_{0}}{z_{0}}\in [0,\overline{w}_{1})$, then 
$$lim_{n\rightarrow\infty} x_{n}=0,$$
$$lim_{n\rightarrow\infty} y_{n}=0,$$
$$lim_{n\rightarrow\infty} z_{n}=\infty.$$
If $\frac{x_{0}+y_{0}}{z_{0}}=\overline{w}_{1}$,
then 
$$ x_{n}=\frac{\overline{w}_{1}}{A_{1}+\overline{w}_{1}},\;\;n\geq 1,$$
$$ y_{n}=\frac{\overline{w}_{1}}{A_{2}+\overline{w}_{1}},\;\;n\geq 1,$$
$$ z_{n}=1+\frac{\alpha }{\overline{w}_{1}},\;\;n\geq 1.$$
If $\frac{x_{0}+y_{0}}{z_{0}}\in (\overline{w}_{1},\infty)$, then 
$$lim_{n\rightarrow\infty} x_{n}=\frac{\overline{w}_{2}}{A_{1}+\overline{w}_{2}},$$
$$lim_{n\rightarrow\infty} y_{n}=\frac{\overline{w}_{2}}{A_{2}+\overline{w}_{2}},$$
$$lim_{n\rightarrow\infty} z_{n}=1+\frac{\alpha }{\overline{w}_{2}}.$$
Whenever $P(w_{m})\leq 0$ or one of the following two conditions hold,
\begin{enumerate}
\item $A_{1}+A_{2}-\alpha A_{1}-\alpha A_{2}-A_{1}A_{2}< 0$ and 
$$2\left(2-A_{1}-A_{2}-\alpha\right)\leq \sqrt{-12\left(A_{1}+A_{2}-\alpha A_{1}-\alpha A_{2}-A_{1}A_{2}\right)},$$
\item $A_{1}+A_{2}-\alpha A_{1}-\alpha A_{2}-A_{1}A_{2}= 0$ and $A_{1}+A_{2}+\alpha \geq 2$,
\end{enumerate}
then 
$$lim_{n\rightarrow\infty} x_{n}=0,$$
$$lim_{n\rightarrow\infty} y_{n}=0,$$
$$lim_{n\rightarrow\infty} z_{n}=\infty.$$
\end{proof}

In our next example, we use the projective change of variables of hyperbolic type to reduce a 3 dimensional system to the rational system in the plane numbered (11,11), in the numbering system introduced in \cite{cklm}. We then further reduce the system to a second order difference equation which decouples into two Riccati difference equations for which the solutions are known. This allows us to determine the qualitative behavior of the original system of 3 rational difference equations.
\begin{expl}
Consider the following system of three rational difference equations,
$$x_{n+1}=\frac{1}{Az_{n}+By_{n}},\;\; n=0,1,2, \dots,$$
$$y_{n+1}=\frac{1}{Cz_{n}+Dx_{n}},\;\; n=0,1,2, \dots,$$
$$z_{n+1}=\frac{1}{z_{n}},\;\; n=0,1,2, \dots,$$
with initial conditions $x_{0},y_{0},z_{0}\in (0,\infty )$ and parameters $A,B,C,D \in (0,\infty )$.
Then every solution converges to the not necessarily prime period 2 solution,
\par\vspace{.2cm}
$$x_{2n}=z_{0}\frac{D-AC-B + \sqrt{(D-AC-B)^{2}+4ADC}}{2AD},$$
\par\vspace{.2cm}
$$y_{2n}=z_{0}\frac{2A}{AC+D-B + \sqrt{(D-AC-B)^{2}+4ADC} },$$
\par\vspace{.2cm}
$$z_{2n}=z_{0},$$
\par\vspace{.2cm}
$$x_{2n+1}=\frac{D-AC-B + \sqrt{(D-AC-B)^{2}+4ADC}}{2ADz_{0}},$$
\par\vspace{.2cm}
$$y_{2n+1}=\frac{2A}{z_{0}\left(AC+D-B + \sqrt{(D-AC-B)^{2}+4ADC}\right) },$$
\par\vspace{.2cm}
$$z_{2n+1}=\frac{1}{z_{0}}.$$
\end{expl}

\begin{proof}

Dividing through by $z_{n+1}$ we have the following relationship,
$$\frac{x_{n+1}}{z_{n+1}}=\frac{1}{A+B\frac{y_{n}}{z_{n}}},\;\; n=0,1,2, \dots,$$
$$\frac{y_{n+1}}{z_{n+1}}=\frac{1}{C+D\frac{x_{n}}{z_{n}}},\;\; n=0,1,2, \dots.$$
Relabeling so that $u_{n}=\frac{x_{n}}{z_{n}}$ and $v_{n}=\frac{y_{n}}{z_{n}}$ we have,
$$u_{n+1}=\frac{1}{A+Bv_{n}},\;\; n=0,1,2, \dots,$$
$$v_{n+1}=\frac{1}{C+Du_{n}},\;\; n=0,1,2, \dots.$$
So the sequence $\{u_{n}\}^{\infty}_{n=1}$ satisfies the following second order rational difference equation,
$$u_{n+1}=\frac{C+Du_{n-1} }{AC+ADu_{n-1}+B},\;\; n=1,2, \dots,$$
with initial conditions $u_{0}=u_{0}$ and $u_{1}=\frac{1}{A+Bv_{0}}$. This equation decouples into two Riccati difference equations. Thus every solution converges to the positive equilibirum 
$$\bar{u}=\frac{D-AC-B + \sqrt{(D-AC-B)^{2}+4ADC}}{2AD}.$$
Thus $\{v_{n}\}^{\infty}_{n=1}$ converges to the positive equilibirum 
$$\bar{v}=\frac{2A}{AC+D-B + \sqrt{(D-AC-B)^{2}+4ADC} }.$$
Now, clearly $z_{2n}=z_{0}$ and $z_{2n+1}=\frac{1}{z_{0}}$ for the sequence $\{z_{n}\}^{\infty}_{n=1}$. Since $\{u_{n}\}^{\infty}_{n=1}$ converges to the positive equilibirum 
$$\bar{u}=\frac{D-AC-B + \sqrt{(D-AC-B)^{2}+4ADC}}{2AD},$$
and $u_{n}=\frac{x_{n}}{z_{n}}$, $\{x_{n}\}^{\infty}_{n=1}$ converges to the period two sequence, 

$$x_{2n}= z_{0}\frac{D-AC-B + \sqrt{(D-AC-B)^{2}+4ADC}}{2AD},$$
and
$$ x_{2n+1}= \frac{D-AC-B + \sqrt{(D-AC-B)^{2}+4ADC}}{2ADz_{0}}.$$
Since $\{v_{n}\}^{\infty}_{n=1}$ converges to the positive equilibirum 
$$\bar{v}=\frac{2A}{AC+D-B + \sqrt{(D-AC-B)^{2}+4ADC} },$$
and $v_{n}=\frac{y_{n}}{z_{n}}$, $\{y_{n}\}^{\infty}_{n=1}$ converges to the period two sequence, 

$$y_{2n}= z_{0}\frac{2A}{AC+D-B + \sqrt{(D-AC-B)^{2}+4ADC} },$$
and
$$y_{2n+1}= \frac{2A}{z_{0}\left(AC+D-B + \sqrt{(D-AC-B)^{2}+4ADC}\right) }.$$
Thus, for our original system, every solution with initial conditions $x_{0}, y_{0}, z_{0}$,
converges to the not necessarily prime period 2 solution,
\par\vspace{.2cm}
$$x_{2n}=z_{0}\frac{D-AC-B + \sqrt{(D-AC-B)^{2}+4ADC}}{2AD},$$
\par\vspace{.2cm}
$$y_{2n}=z_{0}\frac{2A}{AC+D-B + \sqrt{(D-AC-B)^{2}+4ADC} },$$
\par\vspace{.2cm}
$$z_{2n}=z_{0},$$
\par\vspace{.2cm}
$$x_{2n+1}=\frac{D-AC-B + \sqrt{(D-AC-B)^{2}+4ADC}}{2ADz_{0}},$$
\par\vspace{.2cm}
$$y_{2n+1}=\frac{2A}{z_{0}\left(AC+D-B + \sqrt{(D-AC-B)^{2}+4ADC}\right) },$$
\par\vspace{.2cm}
$$z_{2n+1}=\frac{1}{z_{0}}.$$
\end{proof}

\section{Conclusion}
We have presented 3 full families of projective rational difference equations with the corresponding changes of variables. This should be of interest to researchers in the area for two reasons. First, we point out the geometric origins of well known changes of variables. Second, we present new families of projective rational difference equations and give a useful change of variables in each case.  

Moreover, we have presented several special cases where the projective change of variables can be used to completely determine the qualitative behavior. These special cases are interesting in their own right and it would have been exceedingly difficult to determine the qualitative behavior in certain cases without the projective change of variables. 


\begin{thebibliography}{99}
\bibitem{cklm} E. Camouzis, M.R.S. Kulenovi\'c, G. Ladas, and O. Merino, Rational systems in the plane,\;\emph{J. Difference Equ. Appl.}\;\textbf{15}(2009), 303-323.
\vspace{0.1 cm}
\bibitem{riccati} E.A. Grove, G. Ladas, L.C. McGrath, and C.T. Teixeira, Existence and behavior of solutions of a rational system,\;\emph{Commun. Appl. Nonlinear Anal.}\;\textbf{8}(2001), 1-25.
\vspace{0.1 cm}
\bibitem{huangknopf}  Y.S. Huang and P.M. Knopf, Global convergence properties of first-order homogeneous systems of rational difference equations,\;\emph{J. Difference Equ. Appl.}\;, forthcoming article, DOI:10.1080/10236198.2011.590802.
\vspace{0.1 cm}
\bibitem{kulenovicladas} M.R.S. Kulenovi\'c and G. Ladas, \emph{Dynamics of Second Order Rational Difference Equations}, Chapman \& Hall/CRC
Press, Boca Raton, 2002.
\vspace{0.1 cm}
\bibitem{merino}  O. Merino, Global attractivity of the equilibrium of a difference equation: An elementary proof assisted by a computer algebra system,\;\emph{J. Difference Equ. Appl.}\;\textbf{17}(2011), 33-41.
\vspace{0.1 cm}



\end{thebibliography}
\end{document}